\documentclass[10pt,twoside]{article}

\def\YEAR{\year}\newcount\VOL\VOL=\YEAR\advance\VOL by-1995
\def\firstpage{1}\def\lastpage{1000}
\def\received{}\def\revised{}
\def\communicated{}

\makeatletter
\def\magnification{\afterassignment\m@g\count@}
\def\m@g{\mag=\count@\hsize6.5truein\vsize8.9truein\dimen\footins8truein}
\makeatother

\font\eightrm=cmr8
\font\caps=cmcsc10                    
\font\Caps=cmcsc10 scaled \magstep1   

\makeatletter
\@specialpagefalse\headheight=8.5pt
\def\DocMath{}
\renewcommand{\@evenhead}{%
    \ifnum\thepage>\lastpage\rlap{\thepage}\hfill%
    \else\rlap{\thepage}\slshape\leftmark\hfill{\caps\SAuthor}\hfill\fi}%
\renewcommand{\@oddhead}{%
    \ifnum\thepage=\firstpage{\DocMath\hfill\llap{\thepage}}%
    \else{\slshape\rightmark}\hfill{\caps\STitle}\hfill\llap{\thepage}\fi}%
\makeatother

\def\TSkip{\bigskip}
\newbox\TheTitle{\obeylines\gdef\GetTitle #1
\ShortTitle  #2
\SubTitle    #3
\Author      #4
\ShortAuthor #5
\EndTitle
{\setbox\TheTitle=\vbox{\baselineskip=20pt\let\par=\cr\obeylines%
\halign{\centerline{\Caps##}\cr\noalign{\medskip}\cr#1\cr}}%
        \copy\TheTitle\TSkip\TSkip%
\def\next{#2}\ifx\next\empty\gdef\STitle{#1}\else\gdef\STitle{#2}\fi%
\def\next{#3}\ifx\next\empty%
    \else\setbox\TheTitle=\vbox{\baselineskip=20pt\let\par=\cr\obeylines%
    \halign{\centerline{\caps##} #3\cr}}\copy\TheTitle\TSkip\TSkip\fi%
\centerline{\caps #4}\TSkip\TSkip%
\def\next{#5}\ifx\next\empty\gdef\SAuthor{#4}\else\gdef\SAuthor{#5}\fi%
\ifx\received\empty\relax
    \else\centerline{\eightrm Received: \received}\fi%
\ifx\revised\empty\TSkip%
    \else\centerline{\eightrm Revised: \revised}\TSkip\fi%
\ifx\communicated\empty\relax
    \else\centerline{\eightrm Communicated by \communicated}\fi\TSkip\TSkip%
\catcode'015=5}}\def\Title{\obeylines\GetTitle}
\def\Abstract{\begingroup\narrower
    \parskip=\medskipamount\parindent=0pt{\caps Abstract. }}
\def\EndAbstract{\par\endgroup\TSkip}

\long\def\MSC#1\EndMSC{\def\arg{#1}\ifx\arg\empty\relax\else
     {\par\narrower\noindent%
     2000 Mathematics Subject Classification: #1\par}\fi}

\long\def\KEY#1\EndKEY{\def\arg{#1}\ifx\arg\empty\relax\else
        {\par\narrower\noindent Keywords and Phrases: #1\par}\fi\TSkip}

\newbox\TheAdd\def\Addresses{\vfill\copy\TheAdd\vfill
    \ifodd\number\lastpage\vfill\eject\phantom{.}\vfill\eject\fi}
{\obeylines\gdef\GetAddress #1
\Address #2 
\Address #3
\Address #4
\EndAddress
{\def\xs{4.3truecm}\parindent=0pt
\setbox0=\vtop{{\obeylines\hsize=\xs#1\par}}\def\next{#2}
\ifx\next\empty 
     \setbox\TheAdd=\hbox to\hsize{\hfill\copy0\hfill}
\else\setbox1=\vtop{{\obeylines\hsize=\xs#2\par}}\def\next{#3}
\ifx\next\empty 
     \setbox\TheAdd=\hbox to\hsize{\hfill\copy0\hfill\copy1\hfill}
\else\setbox2=\vtop{{\obeylines\hsize=\xs#3\par}}\def\next{#4}
\ifx\next\empty\ 
     \setbox\TheAdd=\vtop{\hbox to\hsize{\hfill\copy0\hfill\copy1\hfill}
                \vskip20pt\hbox to\hsize{\hfill\copy2\hfill}}
\else\setbox3=\vtop{{\obeylines\hsize=\xs#4\par}}
     \setbox\TheAdd=\vtop{\hbox to\hsize{\hfill\copy0\hfill\copy1\hfill}
                \vskip20pt\hbox to\hsize{\hfill\copy2\hfill\copy3\hfill}}
\fi\fi\fi\catcode'015=5}}\gdef\Address{\obeylines\GetAddress}

\begin{document}
\Title   Hasse invariant and group cohomology
\ShortTitle 
\SubTitle   
\Author  Bas Edixhoven, Chandrashekhar Khare
\ShortAuthor 
\EndTitle
\Abstract 

Let $p$ be a prime number. The Hasse invariant is a modular form
modulo~$p$ that is often used to produce congruences between modular
forms of different weights. We show how to produce such congruences
between forms of weights $2$ and $p+1$, in terms of group
cohomology. We also show how our method works in the contexts of
quadratic imaginary fields (where there is no Hasse invariant
available) and Hilbert modular forms over totally real fields of even
degree. 
\EndAbstract
\MSC 
11F33
\EndMSC
\KEY 
\EndKEY
\Address 
Mathematisch Instituut
Universiteit Leiden 
P.O.~Box~9512
2300 RA Leiden 
The Netherlands 
edix@math.leidenuniv.nl
\Address
Department of Mathematics
University of Utah
155 South 1400 East
Salt Lake City
UT~84112
U.S.A.
shekhar@math.utah.edu
\Address
\Address
\EndAddress

\newtheorem{theorem}{Theorem}
\newtheorem{lemma}{Lemma}
\newtheorem{prop}{Proposition}
\newtheorem{cor}{Corollary}
\newtheorem{remark}{Remark}
\newtheorem{example}{Example}
\newtheorem{conjecture}{Conjecture}
\newtheorem{definition}{Definition}
\newtheorem{quest}{Question}

\newcommand{\parab}{{\mathrm{par}}}

\section{The puzzle}
Let $p\geq5$ be a prime number.  In the theory of modular forms
mod~$p$ (see [S] and~[SwD]) a special role is played by the {\it Hasse
invariant} and the $\Theta$~operator.

We fix an embedding $\overline{\bf Q} \hookrightarrow \overline{{\bf
Q}_p}$, and denote the corresponding place of $\overline{\bf Q}$
by~$\wp$.  We have a modular form $E_{p-1}$ of weight $p-1$ in
$M_{p-1}(SL_2({\bf Z}),{\bf Z}_p)$, that is congruent to 1 mod~$\wp$
(see [S] and~[SwD]). By congruence we will mean a congruence of
Fourier coefficients at almost all primes.  The modular form $E_{p-1}$
is the normalised form of the classical Eisenstein series, and has
$q$-expansion
$$
1-2(p-1)/B_{p-1}\sum \sigma_{p-2}(n)q^n,
$$
and the congruence property of $E_{p-1}$ is then a consequence of the
divisibility of the denominator of $B_{p-1}$ by~$p$ (the theorem of
Clausen-von Staudt: see [S, \S1.1]).

Multiplying by $E_{p-1}$ gives the fact that for any positive integer
$N$ prime to~$p$ a weight 2 form in $S_2(\Gamma_1(N),{\bf Z}_p)$ is
congruent mod $p$ to a weight $p+1$ form in $S_{p+1}(\Gamma_1(N),{\bf
Z}_p)$. It follows that a weight 2 normalized eigenform in
$S_2(\Gamma_1(N),\overline{{\bf Z}_p})$ is congruent mod $\wp$ to a
weight $p+1$ eigenform in $S_{p+1}(\Gamma_1(N),\overline{{\bf Z}_p})$
(see~[DS, \S6.10]).

For $N\geq 5$ prime to~$p$, the Hasse invariant constructed
geometrically (see [Kz]) is a global section of the coherent sheaf
$\omega^{\otimes p-1}_{X_1(N)_{{\bf F}_p}}$ where
$\omega_{X_1(N)_{{\bf F}_p}}$ is the pull back of the canonical sheaf
$\Omega_{{\cal E}/X_1(N)_{{\bf F}_p}}$ by the zero section of the map
${\cal E}_{{\bf F}_p}\rightarrow X_1(N)_{{\bf F}_p}$, and with ${\cal
E}_{{\bf F}_p}$ the universal generalised elliptic curve over
$X_1(N)_{{\bf F}_p}$: $E_{p-1}$ can be interpreted as a characteristic
zero lift of the Hasse invariant (Deligne).

The $\Theta$ operator on modular forms mod~$p$ is defined by:
$$
\Theta(\sum a_nq^n)=\sum na_nq^n
$$
where $a_n \in \overline{{\bf F}_p}$. It preserves levels, and
increases weights by $p+1$, i.e., it gives maps:
$$
M_k(\Gamma_1(N),{\bf F}_p) \longrightarrow M_{k+p+1}(\Gamma_1(N),{\bf
F}_p),
$$
preserving cusp forms.  The analog in group cohomology of the $\Theta$
operator on mod $p$ modular forms, can be found in [AS]. The aim of
this note is to find a group theoretic substitute for the Hasse
invariant.

Unlike as is done in [AS] in the case of~$\Theta$, for good reasons we
cannot find an {\it element} in group cohomology that is an analog of
the Hasse invariant.  What we do find instead is a procedure for
raising weights by $p-1$ of mod $p$ Hecke eigenforms of weight two
(preserving the level) that is one of the principal uses of the Hasse
invariant.

Using the Eichler-Shimura isomorphism the relevant Hecke modules are
$H^1(\Gamma_1(N),{\bf F}_p)$ and $H^1(\Gamma_1(N),{\sf
Symm}^{p-1}({\bf F}_p^2))$.  From the viewpoint of group cohomology
the above considerations give that a Hecke system of eigenvalues
$(a_l)_{l\neq p}$ in the former also arises from the latter. This at
first sight is puzzling as indeed the $p-1^{\sf st}$ symmetric power
of the standard 2-dimensional representation of $SL_2({\bf F}_p)$ is
irreducible. In this short note we ``resolve'' this puzzle.

Indeed the solution to the puzzle is implicit in an earlier paper of
one of us (cf. Remark~4 at the end of Section~3 of~[K]) where the
issue arose in trying to understand why the methods for studying
Steinberg lifts of an irreducible modular Galois representation
$\rho\colon G_{\bf Q} \rightarrow GL_2(\overline{{\bf F}_p})$ were
qualitatively different from those for studying principal series and
supercuspidal lifts. The buzzwords there were that the $p$-dimensional
minimal $K$-type of a Steinberg representation of $GL_2({\bf Q}_p)$
also arises in the restriction to $GL_2({\bf Z}_p)$ of any unramified
principal series representation of~$GL_2({\bf Q}_p)$.

The key to the solution of this puzzle (again) is a study of the
degeneracy map $H^1(\Gamma_1(N),{\bf F}_p)^2 \to H^1(\Gamma_1(N) \cap
\Gamma_0(p),{\bf F}_p)$.

In Section~3 we will give applications of our method in the situtation
of imaginary quadratic fields, where the ``geometric Hasse invariant''
perforce is not available. Furthermore modular forms in this setting
do not have a multiplicative structure. We owe this observation, and
indeed the suggestion that our methods should work in this case, to
Ian Kiming. Our cohomological methods do work in this situation, but
have the (inherent) defect that results are about characteristic $p$
modular forms, and may not be used directly to produce congruences
between characteristic $0$ eigenforms of different weights.  This
comes from the fact that in this situation torsion in cohomology can
possibly occur even after localisation at ``interesting'' maximal
ideals of the Hecke algebra (see the concluding remark).

In Section~4 we deal with the case of totally real fields and
in Section~5 we spell out some consequences of~[K] for raising 
of levels in higher weights.

\section{The solution to the puzzle}
Let us recall the hypotheses: $p\geq 5$ is prime, and $N\geq1$ is
prime to~$p$. Consider the cohomology groups $H^1(\Gamma_1(N),{\bf
F}_p)$ and $H^1(\Gamma_1(N) \cap \Gamma_0(p),{\bf F}_p)$.  We have the
standard action of Hecke operators $T_r$ on these cohomology
groups. We recall that we only consider the action for $(r,p)=1$.  We
have the degeneracy map
$$
\alpha:H^1(\Gamma_1(N),{\bf F}_p)^2 \rightarrow 
H^1(\Gamma_1(N) \cap \Gamma_0(p),{\bf F}_p)
$$
that is defined to be the sum $\alpha_1+\alpha_2$ where
$\alpha_1$ is the restiction map, and $\alpha_2$ the ``twisted''
restriction map, given by conjugation by 
$$
g:=\left(\matrix{p&0\cr 0&1\cr}\right)
$$
followed by restriction. The map $\alpha$ is equivariant for the
$T_r$'s that we consider. We have the following variant of a lemma of
Ihara and Ribet (see~[R], and also~[CDT, 6.3.1]).

\begin{lemma}
The map $\alpha\colon H^1(\Gamma_1(N),{\bf F}_p)^2 \rightarrow
H^1(\Gamma_1(N) \cap \Gamma_0(p),{\bf F}_p)$ is injective.
\end{lemma}
\noindent{\bf Proof.} Let $\Delta$ be the subgroup of $SL_2({\bf
Z}[1/p])$ of elements congruent to $({1\atop 0}{*\atop 1})$
modulo~$N$. The arguments of [S2, II,~\S1.4] show that $\Delta$ is the
amalgam of $\Gamma_1(N)$ and $g\Gamma_1(N)g^{-1}$ along their
intersection $\Gamma_1(N)\cap\Gamma_0(p)$. The universal property of
amalgams then implies that the kernel of~$\alpha$ is $H^1(\Delta,{\bf
F}_p)$ i.e., $Hom(\Delta,{\bf F}_p)$. By [S1], each subgroup of finite
index of $SL_2({\bf Z}[1/p])$ is a congruence subgroup, hence each
morphism from $\Delta$ to~${\bf F}_p$ factors through the image
$\Delta_n$ of $\Delta$ in some $SL_2({\bf Z}/n{\bf Z})$ with $n$ prime
to~$p$. The result follows, as $p$ is at least~$5$ and does not
divide~$N$. (We use that $SL_2({\bf Z})$ maps surjectively
to~$SL_2({\bf Z}/n{\bf Z})$.) 

\vspace{3mm}

By Shapiro's lemma we see that $H^1(\Gamma_1(N) \cap \Gamma_0(p),{\bf
F}_p)$ is isomorphic (as a Hecke module) to $H^1(\Gamma_1(N), {\bf
F}_p[{\bf P}^1({\bf F}_p)])$. Using an easy computation of Brauer
characters we deduce that the semisimplification of ${\bf F}_p[{\bf
P}^1({\bf F}_p)]$ under the natural action of $\Gamma_1(N)$ (that
factors through $\Gamma_1(N)/\Gamma_1(N) \cap \Gamma(p)$) is ${\sf id}
\oplus {\sf Symm}^{p-1}({\bf F}_p^2)$. In fact as the cardinality of
${\bf P}^1({\bf F}_p)$ is prime to $p$ we deduce that this is indeed
true even before semisimplification, i.e., ${\bf F}_p[{\bf P}^1({\bf
F}_p)]$ is semisimple as a $SL_2({\bf F}_p)$-module. The submodule
${\sf id}$ is identified with the constant functions, with complement
the functions with zero average. Thus we identify $H^1(\Gamma_1(N)
\cap \Gamma_0(p),{\bf F}_p)$ with $H^1(\Gamma_1(N),{\bf F}_p) \oplus
H^1(\Gamma_1(N),{\sf Symm}^{p-1}({\bf F}_p^2))$. The degeneracy
map~$\alpha$ takes the form:
$$
H^1(\Gamma_1(N),{\bf F}_p)^2 \rightarrow H^1(\Gamma_1(N),{\bf F}_p)
\oplus H^1(\Gamma_1(N),{\sf Symm}^{p-1}({\bf F}_p^2)).
$$

\begin{lemma}
The map:
$$
\beta\colon H^1(\Gamma_1(N),{\bf F}_p) \longrightarrow
H^1(\Gamma_1(N),{\sf Symm}^{p-1}({\bf F}_p^2)),
$$
that is the composition of $\alpha_2$ with the projection of
$H^1(\Gamma_1(N) \cap \Gamma_0(p),{\bf F}_p)$ to $H^1(\Gamma_1(N),{\sf
Symm}^{p-1}({\bf F}_p^2))$, is injective.
\end{lemma}
\noindent{\bf Proof.} This is an immediate consequence of Lemma~1, and
the fact that $\alpha_1\colon H^1(\Gamma_1(N),{\bf F}_p) \rightarrow
H^1(\Gamma_1(N) \cap \Gamma_0(p),{\bf F}_p)$ has image exactly the
first summand of $H^1(\Gamma_1(N),{\bf F}_p) \oplus
H^1(\Gamma_1(N),{\sf Symm}^{p-1}({\bf F}_p^2))$.

\vspace{3mm}

In view of the Eichler-Shimura isomorphism (see [DI, \S12]), we have a
new proof by a purely group cohomological method of the following
well-known result.

\begin{cor}
Suppose moreover that $N\geq 5$. A semi-simple representation
$\rho\colon G_{\bf Q} \rightarrow GL_2(\overline{{\bf F}_p})$ that
arises from $S_2(\Gamma_1(N),\overline{{\bf Q}_p})$ also arises from
$S_{p+1}(\Gamma_1(N),\overline{{\bf Q}_p})$.
\end{cor}

\noindent{\bf Proof.} This follows from the above lemma, together
with the following facts.

\begin{enumerate}

\item The degeneracy map $\beta$ is Hecke equivariant for the $T_r$'s
that we consider, and it sends $H^1_\parab(\Gamma_1(N),{\bf F}_p)$ to
$H^1_\parab(\Gamma_1(N),{\sf Symm}^{p-1}({\bf F}_p^2))$. The subscript
``par'' denotes parabolic cohomology, i.e., the intersection of the
kernels of the restriction maps to the cohomology of the unipotent
subgroups of~$\Gamma_1(N)$.

\item For $V$ any $\Gamma_1(N)$-module that is free of finite rank
over~${\bf Z}$, and such that $H^0(\Gamma_1(N),{\bf F}_p\otimes
V^\vee)=0$, the map:
$$
H^1_\parab(\Gamma_1(N),V)\to 
H^1_\parab(\Gamma_1(N),{\bf F}_p\otimes V)
$$
is surjective (one uses that $H^1_\parab(\Gamma_1(N),{\bf F}_p\otimes
V)$ is a quotient of $H^1_c(Y_1(N),{\bf F}_p\otimes{\cal F}_V)$, with
${\cal F}_V$ the sheaf given corresponding to~$V$, and that
$H^2_c(Y_1(N),{\bf F}_p\otimes{\cal F}_V)=0$ by Poincar\'e duality).
\end{enumerate}

\vspace{3mm}
\noindent{\bf Remarks.}

\begin{enumerate}

\item One can ask the converse question as to which maximal ideals
$\sf m$ of the Hecke algebra acting on $H^1(\Gamma_1(N),{\sf
Symm}^{p-1}({\bf F}_p^2))$ are pull backs of maximal ideals of the
Hecke algebra acting on $H^1(\Gamma_1(N),{\bf F}_p)$. Then, for
non-Eisenstein~$\sf m$, the answer is in terms of the Galois
representation~$\rho_{\sf m}$: a necessary and sufficient condition is
that $\rho_{\sf m}$ be finite flat at~$p$ (see~[R1, Thm.~3.1]).
Perhaps one does not expect to have a group cohomological approach to
such a subtle phenomenon.

\item Let $\rho$ be an irreducible 2-dimensional mod $p$
representation of~$G_{\bf Q}$. We just saw that if $\rho$ arises from
$S_2(\Gamma_1(N))$, then it also does so from $S_{p+1}(\Gamma_1(N))$,
using only group cohomology. See~[He] for the case of higher weights,
when the raising of weights by $p-1$ is explained more directly in
terms of the Jordan-H\"older series of the ${\sf Symm}^n({\bf
F}_p^2)$. On the other hand the method here being cohomological cannot 
be used to raise weights from 1 to~$p$, while multiplication by the Hasse
invariant can be used for this (see~[DS]).
\end{enumerate}

\section{Imaginary quadratic fields}
Let $p\geq 5$ be a prime number, $K$ an imaginary quadratic field in
which $p$ is inert, and $\cal N$ a non-zero ideal in the ring
of integers ${\cal O}_K$ not containing~$p$.

Let $\Gamma_1({\cal N})$ be the congruence subgroup of $SL_2({\cal
O}_K)$ of level~${\cal N}$. As before we have the degeneracy map:
$$
\alpha\colon H^1(\Gamma_1({\cal N}),{\bf F}_p)^2 \rightarrow 
H^1(\Gamma_1({\cal N}) \cap \Gamma_0(p),{\bf F}_p)
$$
that is defined to be the sum $\alpha_1 \oplus \alpha_2$ where
$\alpha_1$ is the restriction map, and $\alpha_2$ the ``twisted''
restriction map, given by ``conjugation'' by 
$$
g:=\left(\matrix{p&0\cr 0&1\cr}\right)
$$
followed by restriction.  Then we again have:

\begin{lemma}
The map:
$\alpha\colon H^1(\Gamma_1({\cal N}),{\bf F}_p)^2 \rightarrow
H^1(\Gamma_1({\cal N}) \cap \Gamma_0(p),{\bf F}_p)$
is injective.
\end{lemma}
\noindent{\bf Proof.} One replaces ${\bf Z}$ by ${\cal O}_K$ and $N$
by~${\cal N}$ in the proof of Lemma~1. Strong approximation (see~[PR])
guarantees that the reduction map from $SL_2({\cal O}_K)$ to
$SL_2({\cal O}_K/n{\cal O}_K)$ is surjective for all $n\geq1$.

\vspace{3mm}

By Shapiro's lemma we see that $H^1(\Gamma_1({\cal
N})\cap\Gamma_0(p),{\bf F}_p)$ is isomorphic to $H^1(\Gamma_1({\cal
N}), {\bf F}_p[{\bf P}^1({\bf F}_{\wp})])$, where ${\bf F}_{\wp}={\cal
O}_K/\wp$. Using an easy computation of Brauer characters we deduce
that the semisimplification of ${\bf F}_\wp[{\bf P}^1({\bf F}_{\wp})]$
under the natural action of $\Gamma_1({\cal N})$ (that factors through
$\Gamma_1({\cal N})/ \Gamma_1({\cal N}) \cap \Gamma(\wp)$) is ${\sf
id} \oplus {\sf Symm}^{p-1}({\bf F}_\wp^2) \otimes {\sf
Symm}^{p-1}({\bf F}_\wp^2)^{\sigma}$, with $\sigma$ the non-trivial
automorphism of~${\bf F}_{\wp}$, and the superscript denotes that the
action has been twisted by~$\sigma$. Note that ${\sf Symm}^{p-1}({\bf
F}_\wp^2) \otimes {\sf Symm}^{p-1}({\bf F}_\wp^2)^{\sigma}$ is
irreducible as a $SL_2({\bf F}_{\wp})$-module: this is a particular
case of the well-known tensor product theorem of Steinberg (see~[St]).
In fact as the cardinality of ${\bf P}^1({\bf F}_{\wp})$ is prime to
$p$ we deduce as before that this is indeed true even before
semisimplification, i.e., ${\bf F}_p[{\bf P}^1({\bf F}_{\wp})]$ is
semisimple as a $SL_2({\bf F}_{\wp})$-module.

Thus $\alpha$ maps $H^1(\Gamma_1({\cal N}),{\bf F}_\wp)^2$ into the
direct sum of $H^1(\Gamma_1({\cal N}),{\bf F}_\wp)$ and
$H^1(\Gamma_1({\cal N}),{\sf Symm}^{p-1}({\bf F}_\wp^2)\otimes {\sf
Symm}^{p-1}({\bf F}_\wp^2)^{\sigma})$, and composing with the
projection to the second term gives a map:
$$
\beta\colon H^1(\Gamma_1({\cal N}),{\bf F}_\wp) \rightarrow
H^1(\Gamma_1({\cal N}),{\sf Symm}^{p-1}({\bf F}_\wp^2) \otimes
{\sf Symm}^{p-1}({\bf F}_\wp^2)^{\sigma})
$$ 

\begin{lemma}
The map $\beta$ is injective.
\end{lemma}

\noindent{\bf Proof.} After the above discussion, this is an immediate
consequence of Lemma~3 as before.

\vspace{3mm}

The map~$\beta$ is equivariant for the the action of all Hecke
operators outside $p$ (i.e., induced by elements of $SL_2({\bf
Q}_l\otimes K)$ for $l\neq p$). Thus we have proved:
\begin{cor}
Each system of Hecke eigenvalues in $\overline{{\bf F}_p}$ that arises from
$H^1(\Gamma_1({\cal N}),{\bf F}_p)$ also arises from
$H^1(\Gamma_1({\cal N}),{\sf Symm}^{p-1}({\bf F}_\wp^2) \otimes {\sf
Symm}^{p-1}({\bf F}_\wp^2)^{\sigma})$.
\end{cor}

\vspace{3mm}

\noindent{\bf Remarks.} 
\begin{enumerate}
\item This result as it stands does not yield any
information about congruences of systems of Hecke eigenvalues
occurring in characteristic zero, as there is a problem with lifting.
More precisely, the obstruction is in the $p$-torsion of
$H^2(\Gamma_1({\cal N}),{\sf Symm}^{p-1}({\cal O}^2) \otimes {\sf
Symm}^{p-1}({\cal O}^2)^{\sigma})$.
\item The condition that $p$ be split in~${\cal O}_K$ is probably not
essential. 
\end{enumerate}

\section{Totally real fields}

The method is also applicable in the case of Hilbert modular forms for
totally real fields of even degree. We quickly sketch the approach
which is similar to that of the previous two sections. Let $F/{\bf Q}$
be a totally real, cyclic extension of even degree, ${\rm Gal}(F/{\bf
Q})=\langle \sigma \rangle$, $p \geq 5$ an inert prime, with $\wp$ the
unique prime of $F$ above it, ${\bf F}_{\wp}$ the residue field
at~$\wp$, and $\cal N$ an ideal of the ring of integers of $F$ that is
prime to~$p$.

Consider the quaternion algebra $D$ over $F$ ramified at all infinite
places and unramified at all finite places, and for any
$F$-algebra~$R$, set $B(R)=(D \otimes_F R)^*$. Let ${\bf A}$ be the
adeles of~$F$, and $U_1({\cal N})$ the standard open compact (mod
centre) subgroup of~$B({\bf A})$. The space of mod $p$ weight 2
modular forms $S({\cal N})$ (resp., $S({\cal N},\wp)$) for $U_1({\cal
N})$ (resp., $U_1({\cal N}) \cap U_0(\wp)$) in this case consists of
functions $B({\bf A}) \rightarrow {\bf F}_p$ that are left and right
invariant under $B(F)$ and $U_1({\cal N})$ (resp., $U_1({\cal N}) \cap
U_0(\wp)$) respectively, modulo the space of functions that factor
through the norm. These spaces come equipped with Hecke actions. This
time controlling the kernel of the degeneracy map $S({\cal N})^2
\rightarrow S({\cal N},\wp)$, i.e., analog of Lemmas~1 and~3, is
easier and follows from strong approximation.  Note again that the
representation ${\sf Symm}^{p-1}({\bf F}_\wp^2) \otimes {\sf
Symm}^{p-1}({\bf F}_\wp^2)^{\sigma} \otimes \cdots \otimes {\sf
Symm}^{p-1}({\bf F}_\wp^2)^{\sigma^{[F:{\bf Q}]-1}}$ of $GL_2({\bf
F}_{\wp})$ (which again is a direct summand, with complement the
trivial representation, of the induction of the trivial representation
from the Borel subgroup of $GL_2({\bf F}_{\wp})$ to $GL_2({\bf
F}_{\wp})$) is irreducible as a consequence of Steinberg's tensor
product theorem.

Now following the method of the previous section, and invoking the
Jacquet-Langlands correspondence yields the following result.

\begin{prop}
With notation as above, suppose that an irreducible representation
$\rho\colon G_F \rightarrow GL_2(\overline{{\bf F}_p})$ arises from a
Hilbert modular form on $\Gamma_1({\cal N})$ of weight $(2,\ldots,2)$.
Then it also arises from a a Hilbert modular form on $\Gamma_1({\cal
N})$ of weight $(p+1,\ldots,p+1)$.
\end{prop}

\noindent{\bf Remarks:} It will be of interest to work out some of the
Hasse invariants with non-parallel weights obtained by E.~Goren from
the viewpoint of this paper (see~[G]).

\section{Congruences between forms of level $N$ and $Np$ for weights
$2 < k \leq p+1$}

We take this opportunity to write down a level raising criterion from
level $N$ to level $Np$ for all weights that is easily deduced from
Corollary~9 of~[K], and list some errata to~[K].

\begin{prop}
Let $f$ be a newform in $S_k(\Gamma_1(N))$ for an integer $k \geq 2$, such
that the mod $\wp$ representation corresponding to it is irreducible. Then:
\begin{itemize}
\item If $k=2$, $f$ is congruent to a $p$-new form in $S_k(\Gamma_1(N)
\cap \Gamma_0(p))$ if and only if $a_p(f)^2=\varepsilon_f(p)$ mod
$\wp$ where $\varepsilon_f$ is the nebentypus of $f$.
\item If $2< k \leq p+1$, $f$ is congruent to a $p$-new form in
$S_k(\Gamma_1(N) \cap \Gamma_0(p))$ if and only if $a_p(f)$ is 0 mod
$\wp$.
\item If $k>p+1$, $f$ is always congruent to a $p$-new form in
$S_k(\Gamma_1(N) \cap \Gamma_0(p))$.
\end{itemize}
\end{prop}

\subsection*{Errata to [K]}

One of us (C.K.) would like to point out some typos in [K]:

\begin{enumerate}
\item
Lines 12 and 19 of Definition 10 page 143 of [K] replace
$$
f\colon D({\bf Q})\backslash D({\bf A}^{\infty})/V
\longrightarrow Hom_{\cal O}(M,Symm^{k-2}({\cal O})).
$$
by
$$
f\colon D({\bf Q})\backslash D({\bf A}^{\infty})
\longrightarrow Hom_{\cal O}(M,Symm^{k-2}({\cal O})).
$$
\item On line 12, page 146 of [K] replace $V_1(N)^p$ by $V_1(N)^p
\times V_p$.
\end{enumerate}

\section{Acknowledgements} The second author thanks 
Ian Kiming for helpful conversations during his visit to the
University of Copenhagen. He thanks CNRS for support during his visit
to IRMAR, Universit\'e de Rennes~I.

\section{References}

\noindent [AS] A.~Ash, G.~Stevens, {\it Modular forms in
characteristic $l$ and special values of their $L$-functions}, Duke
Math. J. 53 (1986), no. 3, 849--86

\vspace{3mm}

\noindent [CDT] B.~Conrad, F.~Diamond, R.~Taylor. {\it Modularity of
certain potentially Barsotti-Tate Galois representations.}
J.A.M.S.~{\bf 12} (1999), 521--567.

\vspace{3mm}

\noindent [DI] F.~Diamond, J.~Im {\it Modular forms and modular
curves.} In ``Seminar on Fermat's last theorem'', Canadian
Mathematical Society Conference Proceedings 17, 1995 (V.\ Kumar Murty,
editor).

\vspace{3mm}

\noindent [DS] P.\ Deligne and J-P.\ Serre. {\it Formes modulaires de
poids~1.} Ann.\ Sci.\ Ecole Norm.\ Sup.~(4) 7, 507--530 (1974).

\vspace{3mm}

\noindent [G] E.Z.\ Goren. {\it Hasse invariants for Hilbert modular
varieties.} Israel J.\ Math.~122 (2001), 157--174. 

\vspace{3mm}

\noindent [He] A.\ Herremans. {\it A combinatorial interpretation of
Serre's conjecture on modular Galois representations}, Preprint
2002-03, Orsay. 

\vspace{3mm}

\noindent [K] C.~Khare, {\it A local analysis of congruences in the
$(p,p)$ case: Part II}, Invent. Math. 143 (2001), 129--155.

\vspace{3mm}

\noindent [Kz] N.~Katz, {\it A result on modular forms in
characteristic $p$}, Modular functions of one variable V, pp. 53--61,
Lecture Notes in Math., Vol. 601, Springer, Berlin, 1977.

\vspace{3mm}

\noindent [PR] V.\ Platonov and A.\ Rapinchuk. {\it Algebraic groups
and number theory}. Academic Press, 1994.

\vspace{3mm}

\noindent [R] K.A.~Ribet, {\it Congruence relations between modular
forms}, Proceedings of the International Congress of Mathematicians,
Vol. 1, 2 (Warsaw, 1983), 503--514, PWN, Warsaw, 1984.

\vspace{3mm}

\noindent [R1] K.A.~Ribet, {\it Report on mod $l$ representations of
${\rm Gal}(\overline{\bf Q}/{\bf Q})$.} Proceedings of Symposia in
Pure Mathematics, {\bf 55} (1994), Part~2.

\vspace{3mm}

\noindent [S] J-P.~Serre, {\it Formes modulaires et fonctions zeta
$p$-adiques}, Modular functions of one variable III, pp. 191--268,
Lecture Notes in Math., Vol. 350, Springer, 1973.

\vspace{3mm}

\noindent [S1] J-P.~Serre, {\it Le probleme des groupes de congruence
pour $SL_2$}, Annals of Math. 92 (1970), 489--527.

\vspace{3mm}

\noindent [S2] J-P.~Serre, {\it Trees.} Springer-Verlag, 1980.

\vspace{3mm}

\noindent [St] R.\ Steinberg, {\it Tensor product theorems.}
Proceedings of Symposia in Pure Mathematics Volume~47 (1987),
331--338. 

\vspace{3mm}

\noindent [SwD] H.P.F.~Swinnerton-Dyer, {\it On $\ell$-adic
representations and congruences for coefficients of modular forms},
Modular functions of one variable III, pp. 1--55, Lecture Notes in
Math., Vol. 350, Springer, 1973.

\Addresses
\end{document}